\documentclass{article}

\usepackage{amsthm,amsmath,amssymb}
\usepackage{tikz,url}
\usepackage{mathpazo}
\usepackage{fancybox}

\newcounter{row2}
\newcounter{row5}
\newcounter{row6}
\newcounter{row7}
\newcounter{col2}
\newcounter{col7}
\newcounter{col6}
\newcounter{col5}

\newcommand\setrowFIVE[5]{
  \setcounter{col5}{1}
  \foreach \n in {#1, #2, #3, #4, #5} {
    \edef\x{\value{col5} - 0.5}
    \edef\y{5.5 - \value{row5}}
    \node[anchor=center] at (\x, \y) {\n};
    \stepcounter{col5}
  }
  \stepcounter{row5}
}

\newcommand\setrowTWO[2]{
  \setcounter{col2}{1}
  \foreach \n in {#1, #2} {
    \edef\x{\value{col2} - 0.5}
    \edef\y{2.5 - \value{row2}}
    \node[anchor=center] at (\x, \y) {\n};
    \stepcounter{col2}
  }
  \stepcounter{row2}
}

\newcommand\setrowSIX[6]{
  \setcounter{col6}{1}
  \foreach \n in {#1, #2, #3, #4, #5, #6} {
    \edef\x{\value{col6} - 0.5}
    \edef\y{6.5 - \value{row6}}
    \node[anchor=center] at (\x, \y) {\n};
    \stepcounter{col6}
  }
  \stepcounter{row6}
}

\newcommand\setrowSEVEN[7]{
  \setcounter{col7}{1}
  \foreach \n in {#1, #2, #3, #4, #5, #6, #7} {
    \edef\x{\value{col7} - 0.5}
    \edef\y{7.5 - \value{row7}}
    \node[anchor=center] at (\x, \y) {\n};
    \stepcounter{col7}
  }
  \stepcounter{row7}
}

\newcommand{\circl}{\hspace{4.3mm}\vspace{3mm}\circle{12}}
\newtheorem{theorem}{Theorem}

\begin{document}
\title{The Pulsar Sequence}
\author{Vadim Ponomarenko$^\star$}
\date{$^\star$San Diego State University, \url{vponomarenko@sdsu.edu}\\[2ex] \today }

\maketitle

A recent video \cite{ctc} on the popular YouTube sudoku channel ``Cracking The Cryptic'' featured a non-sudoku puzzle, by the constructor Pulsar.  The task was, subject to a particular circle restriction, to construct a $9\times 9$ Latin square.  A Latin square of size $n$ is a rectangular grid filled with the numbers $1$ to $n$, so that each row and column contains distinct numbers.  Some of the squares (in a spiral pattern) are marked with circles, and the restriction for this puzzle is that a digit in a circle indicates the number of circles that contain that digit.  Perhaps surprisingly, there is a unique solution.

It turns out that this is just one of infinitely many such puzzles, because a similar puzzle has a unique solution, for every size (not just $9\times 9$).  Further, these puzzles are all related to the  mysterious sequence of positive integers
\[1,2,1,3,2,1,4,2,3,1,5,2,3,4,1,6,2,4,3,5,1,7,2,5,4,3,6,1,8,2,6,5,6,3,\ldots\]
We call this the Pulsar sequence, and we investigate its connection to the Pulsar puzzles below.  Readers wishing to try the puzzles without spoilers should pause and solve the samples below before continuing.\\

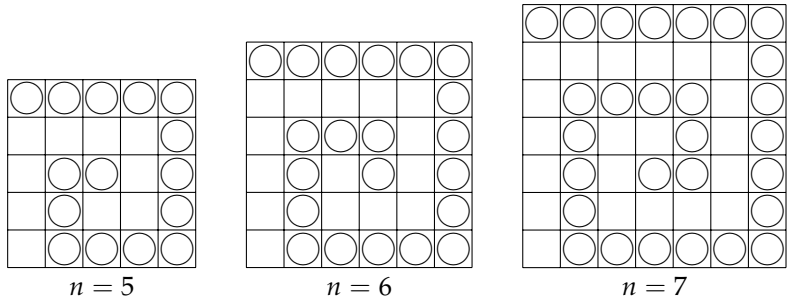
\begin{tikzpicture}[scale=.5]
  \begin{scope}
    \draw (0, 0) grid (5, 5);

    \setcounter{row5}{1}
    \setrowFIVE {\circl }{ \circl}{\circl }  {\circl }{\circl}
    \setrowFIVE {}{}{}{}{\circl}
    \setrowFIVE {}{\circl}{\circl}{}{\circl}
    \setrowFIVE {}{\circl}{}{}{\circl}
    \setrowFIVE {}{\circl}{\circl}{\circl}{\circl}
    \node[anchor=center] at (2.5, -0.5) {$n=5$};
  \end{scope} 
    
\end{tikzpicture} \begin{tikzpicture}[scale=.5]

  \begin{scope}
    \draw (0, 0) grid (6, 6);

    \setcounter{row6}{1}
    \setrowSIX {\circl }{ \circl}{\circl }  {\circl }{\circl}{\circl}
    \setrowSIX {}{}{}{}{}{\circl}
    \setrowSIX {}{\circl}{\circl}{\circl}{}{\circl}
    \setrowSIX {}{\circl}{}{\circl}{}{\circl}
    \setrowSIX {}{\circl}{}{}{}{\circl}
    \setrowSIX {}{\circl}{\circl}{\circl}{\circl}{\circl}

    \node[anchor=center] at (3, -0.5) {$n=6$};
  \end{scope}
\end{tikzpicture} \begin{tikzpicture}[scale=.5]

  \begin{scope}
    \draw (0, 0) grid (7, 7);

    \setcounter{row7}{1}
    \setrowSEVEN{\circl}{\circl}{\circl }  {\circl }{\circl}{\circl}{\circl}
    \setrowSEVEN{}{}{}{}{}{}{\circl}
    \setrowSEVEN{}{\circl}{\circl}{\circl}{\circl}{}{\circl}
    \setrowSEVEN{}{\circl}{}{}{\circl}{}\circl{}
    \setrowSEVEN{}{\circl}{}{\circl}{\circl}{}{\circl}
    \setrowSEVEN{}{\circl}{}{}{}{}{\circl}
    \setrowSEVEN{}{\circl}{\circl}{\circl}{\circl}{\circl}{\circl}

    \node[anchor=center] at (3.5, -0.5) {$n=7$};
  \end{scope}
\end{tikzpicture}

We now make precise which cells are in the spiral, for a given $n$.  There will be $n$ pieces to the spiral, of size $n, n-1, \ldots, 1$, for a total of $\frac{n(n+1)}{2}$ cells with circles.  Each piece will be consecutive cells in a single row or column.  The first piece will be the first row of the grid (size $n$).  The second piece will be the entire $n$-th column of the grid, apart from the top square.  The third piece will be the entire $n$-th row of the grid, apart from the first and last square, and so on.  We mark the pieces with blue lines in the examples below.  We also mark the uncircled squares with red lines; they can be seen to also form a spiral.\\

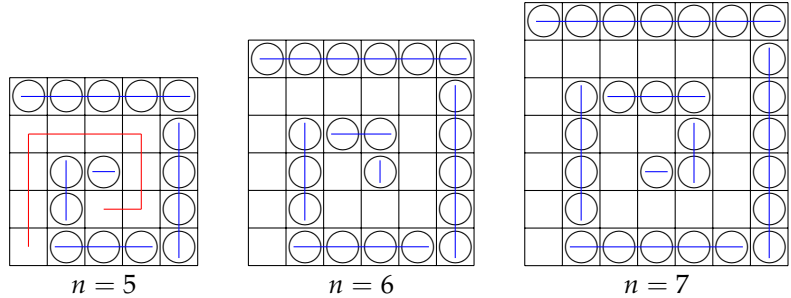
\begin{tikzpicture}[scale=.5]
  \begin{scope}
    \draw (0, 0) grid (5, 5);

    \setcounter{row5}{1}
    \setrowFIVE {\circl }{ \circl}{\circl }  {\circl }{\circl}
    \setrowFIVE {}{}{}{}{\circl}
    \setrowFIVE {}{\circl}{\circl}{}{\circl}
    \setrowFIVE {}{\circl}{}{}{\circl}
    \setrowFIVE {}{\circl}{\circl}{\circl}{\circl}
    \node[anchor=center] at (2.5, -0.5) {$n=5$};
      \draw [blue](0.3,4.5) -- (4.8,4.5);
      \draw [blue](4.5,3.8) -- (4.5,0.2);
      \draw [blue](1.2,0.5) -- (3.8,0.5);
      \draw [blue](1.5,1.2) -- (1.5,2.8);
      \draw [blue](2.2,2.5) -- (2.8,2.5);
      \draw [red](0.5,0.5) -- (0.5,3.5);
      \draw[red](0.5,3.5)--(3.5,3.5);
      \draw[red](3.5,3.5)--(3.5,1.5);
      \draw[red](2.5,1.5)--(3.5,1.5);
    
  \end{scope}
\end{tikzpicture} \begin{tikzpicture}[scale=.5]

  \begin{scope}
    \draw (0, 0) grid (6, 6);

    \setcounter{row6}{1}
    \setrowSIX {\circl }{ \circl}{\circl }  {\circl }{\circl}{\circl}
    \setrowSIX {}{}{}{}{}{\circl}
    \setrowSIX {}{\circl}{\circl}{\circl}{}{\circl}
    \setrowSIX {}{\circl}{}{\circl}{}{\circl}
    \setrowSIX {}{\circl}{}{}{}{\circl}
    \setrowSIX {}{\circl}{\circl}{\circl}{\circl}{\circl}

    \node[anchor=center] at (3, -0.5) {$n=6$};
     \draw [blue](0.3,5.5) -- (5.8,5.5);
          \draw [blue](5.5,4.8) -- (5.5,0.2);
          \draw [blue](1.2,0.5) -- (4.8,0.5);
          \draw [blue](1.5,1.2) -- (1.5,3.8);
          \draw [blue](2.2,3.5) -- (3.8,3.5);
          \draw [blue](3.5,2.2) -- (3.5,2.8);
  \end{scope}
\end{tikzpicture} \begin{tikzpicture}[scale=.5]

  \begin{scope}
    \draw (0, 0) grid (7, 7);

    \setcounter{row7}{1}
    \setrowSEVEN{\circl}{\circl}{\circl }  {\circl }{\circl}{\circl}{\circl}
    \setrowSEVEN{}{}{}{}{}{}{\circl}
    \setrowSEVEN{}{\circl}{\circl}{\circl}{\circl}{}{\circl}
    \setrowSEVEN{}{\circl}{}{}{\circl}{}\circl{}
    \setrowSEVEN{}{\circl}{}{\circl}{\circl}{}{\circl}
    \setrowSEVEN{}{\circl}{}{}{}{}{\circl}
    \setrowSEVEN{}{\circl}{\circl}{\circl}{\circl}{\circl}{\circl}

    \node[anchor=center] at (3.5, -0.5) {$n=7$};
    		\draw [blue](0.3,6.5) -- (6.8,6.5);
              \draw [blue](6.5,5.8) -- (6.5,0.2);
              \draw [blue](1.2,0.5) -- (5.8,0.5);
              \draw [blue](1.5,1.2) -- (1.5,4.8);
              \draw [blue](2.2,4.5) -- (4.8,4.5);
              \draw [blue](4.5,2.2) -- (4.5,3.8);
              \draw[blue] (3.2,2.5) -- (3.8,2.5);

  \end{scope}
\end{tikzpicture}

Note first that the puzzle rules require us to fill these circles with exactly one $1$, two $2$'s, \ldots $n$ $n$'s.  Since each piece is on a row or column of a Latin square, it can contain at most one copy of any number.  Since there are $n$ pieces that must contain $n$ copies of $n$, each piece must contain exactly one $n$.  Now the piece of size $1$ is full, and all remaining $n-1$ pieces must contain exactly one $n-1$.  Continuing in this way, we see that the piece of size $i$ must contain the $i$ numbers $n, n-1, \ldots, n-(i-1)$, in some order.

It will be helpful to define the dual function $(x)_n$ on $\{1,2,\ldots,n\}$ via $(x)_n=n+1-x$.  Note that this reverses the usual order on this set, and applying the dual twice returns what we started with, i.e. $((x)_n)_n=x$.  Importantly, taking the dual of each entry of a Latin square yields another Latin square, which we call the dual Latin square.  Of course, the dual Latin square from a Pulsar puzzle solution no longer satisfies the circle condition, but never mind.  The benefit of considering the dual to a Pulsar puzzle solution, is that in the dual the piece of size $i$ must contain the $i$ numbers $1,2,\ldots, i$, in some order.  

We will show that it is not just ``some order'', but always the same order, and if we append these numbers (in the dual to a Pulsar puzzle solution) from the smallest piece to the largest, we will get precisely terms from the start of the Pulsar sequence.  This sequence comes in blocks of size $1,2,3,\ldots$, as follows:
\[\underbrace{1},\underbrace{2,1},\underbrace{3,2,1},\underbrace{4,2,3,1},\underbrace{5,2,3,4,1},\underbrace{6,2,4,3,5,1},\underbrace{7,2,5,4,3,6,1},\underbrace{8,2,6,5,6,3,7,1},\ldots\]

The $i$'th block contains $a_1, a_2, \ldots, a_i$ and satisfies the symmetric sum property $a_1+a_i=a_2+a_{i-1}=\cdots = i+1$, i.e. the third block satisfies $3+1=2+2=4$ and the fourth block satisfies $4+1=2+3=5$.

For example, consider $n=5$. The dual Pulsar sequence is $(1)_5=5, (2)_5=4,\ldots$, so we can  fill in the circled spiral with   as $5,4,5,3,4,5,2,4,3,5,1,4,3,2,5$.  We fill in the uncircled spiral (in red above)  with $1,2,1,3,2,1,4,2,3,1,5,2,3,4,1$. Hence the full solution is a Pulsar sequence and a dual Pulsar sequence, intertwined as two spirals.

\begin{theorem}
The $n\times n$ Pulsar puzzle has a unique solution, where the circled spiral is filled with the dual of the Pulsar sequence (starting from the center) and the uncircled spiral  is filled with the Pulsar sequence (starting from the center).  Further, in each row and column, each circled number is greater than each uncircled number.
\end{theorem}
\noindent \emph{Proof.} 

\!\!\!\!\!\!\!\!\!\!\!\!\!\begin{minipage}{\textwidth}
\begin{tabular}{ll}
\parbox[l]{4in}{
The proof is by induction on $n$.  As base case we will take $n=2$ (to avoid worrying about trivialities), whose unique solution is pictured at right. Checking all the desired conditions is routine.\\[-8pt]} &
\parbox[t]{0.6in}{
\begin{minipage}{1in}
\begin{tikzpicture}[scale=.5]
    \begin{scope}
      \draw (0, 0) grid (2, 2);
      \setcounter{row2}{1}
      \setrowTWO {\circl }{ \circl}
      \setrowTWO {}{\circl }
      \node[anchor=center] at (1, 0) {};
    \end{scope} 
     \begin{scope}
        \draw (0, 0) grid (2, 2);
        \setcounter{row2}{1}
        \setrowTWO {2}{1}
        \setrowTWO {1}{2}
        \node[anchor=center] at (1,0) {};
      \end{scope}       
  \end{tikzpicture}\end{minipage}}  
  \end{tabular}\end{minipage}
  
Now take $n>2$ and consider its Pulsar puzzle $P$.  The first row consists entirely of circles, while the remainder of the first column is entirely non-circles.  Removing those gives us an $(n-1)\times (n-1)$ Pulsar puzzle  $P'$, which is itself a Pulsar puzzle rotated $90^\circ$ clockwise.  By the inductive hypothesis, there is a unique solution where the circle spiral is filled with the $(n-1)$-dual of the Pulsar sequence.

Removing instead $P$'s first and last rows, and swapping circled cells for non-circled cells, we get $P''$, which is also an $(n-1)\times (n-1)$ Pulsar puzzle, rotated $90^\circ$ counterclockwise.  By the inductive hypothesis, its rows and columns also contain distinct numbers of circles, and  there is a unique solution where the circle spiral is filled with the $(n-1)$-dual of the Pulsar sequence.  We illustrate $P$, $P'$, $P''$ below.

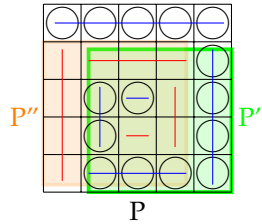
\begin{center}\begin{tikzpicture}[scale=.5]
  \begin{scope}
 \filldraw[fill=green!15, draw=green!100, line width=0.5mm, draw opacity=80](1.2,0) rectangle (5,3.8);
 \filldraw[fill=orange!55, draw=orange!100, line width=0.5mm, nearly transparent](0,0.2) rectangle (3.8,4);

        \draw (0, 0) grid (5, 5);
   \node[color=green] at (5.5,2) {P'};
   \node[color=orange] at (-0.5,2) {P''};
    \setcounter{row5}{1}
    \setrowFIVE {\circl }{ \circl}{\circl }  {\circl }{\circl}
    \setrowFIVE {}{}{}{}{\circl}
    \setrowFIVE {}{\circl}{\circl}{}{\circl}
    \setrowFIVE {}{\circl}{}{}{\circl}
    \setrowFIVE {}{\circl}{\circl}{\circl}{\circl}
    \node[anchor=center] at (2.5, -0.5) {P};
      \draw [blue](0.3,4.5) -- (4.8,4.5);
      \draw [blue](4.5,3.8) -- (4.5,0.2);
      \draw [blue](1.2,0.5) -- (3.8,0.5);
      \draw [blue](1.5,1.2) -- (1.5,2.8);
      \draw [blue](2.2,2.5) -- (2.8,2.5);
      \draw [red](0.5,0.3) -- (0.5,3.8);
      \draw[red](1.2,3.5)--(3.8,3.5);
      \draw[red](3.5,2.8)--(3.5,1.2);
      \draw[red](2.2,1.5)--(2.8,1.5);

  \end{scope}
\end{tikzpicture}\end{center}
  

Next, we produce a solution for $P$.  The solution for $P'$ has one circled $1$, two circled $2$'s, \ldots, $n-1$ circled $(n-1)$'s. The $i$'th circled term is $(x_i)_{n-1}$, where $x_i$ is the $i$'th term of the Pulsar sequence.  What we want instead is $(x_i)_n=(((x_i)_{n-1})_{n-1})_n=(x_i)_{n-1}+1$, so we add one to the circled elements of $P'$ to find the circled elements of $P$.  Now we have the top row of $P$ with empty circles, and also one circled $2$, two circled $3$'s, \ldots, $n-1$ circled $n$'s. 

Comparing $P''$ with $P$, we find that the uncircled cells of $P$ contain the Pulsar sequence, i.e. $n-1$ pieces that contain $n-1$ uncircled $1$'s, $n-2$ uncircled $2$'s, \ldots, $1$ uncircled $n-1$.  We fill those uncircled numbers into $P$.

We have filled in all of $P$ but the top row, and included $n-1$ of each number. We now prove this is a partial Latin square. $P'$ was a Latin square, and we added 1 to each circled number.  However, by our inductive hypothesis all circled values were greater than all uncircled values in each row and column, so adding 1 to the former prevents any duplicate value in any row or column in the portion of $P$ drawn from $P'$.  A similar argument proves that there are no duplicate values in the portion of $P$ drawn from $P''$.  

To prove we have a partial Latin square, we will now prove that each entry in the first column is smaller than each entry in the last column.  Consider the $j$'th row (for $2\le j\le n$).  The first column is the $(n-1)$-th block of the Pulsar sequence, so the corresponding element in $P$ is $a_{j-1}$.  The last column is the $n$-dual of the $(n-1)$-th block of the Pulsar sequence, but reversed.  Hence the corresponding element in $P$ is $n+1-a_{n+1-j}$.  Subtracting, we get $n+1-a_{n+1-j}-a_{j-1}=n+1-n=1$, so in fact each entry in the first column is one smaller than each entry in the last column (in the part filled in, rows $2$ through $n$).

 Hence, considering each column individually, there is a unique way to fill in the top row to make a Latin square.  In fact, those values are all one higher than the number of uncircled squares in their columns, since those are the values made available by adding 1 to the circled entries in $P'$. We see that $P$ contains one circled $1$ (in the top row), two circled $2$'s (one in the top row, one in $P'$), \ldots, $n$ circled $n$'s, so $P$ also satisfies the circle restriction.  
 
 Lastly, we verify the symmetric sum property on the top row, the dual to the new block of the Pulsar sequence.  By symmetry, the number of circled entries in column $i$ is the number of uncircled entries in column $n+1-i$.  Hence the top row satisfies $b_i+b_{n-i}=n+1$.  Since these are duals of the Pulsar sequence, we get $n+1=(a_i)_n+(a_{n+1-i})_n=(n+1-a_i)+(n+1-a_{n+1-i})$, which rearranges to $a_i+a_{n+1-i}=n+1$, as desired.   \qed

~\\[-10pt]


The Pulsar sequence has some clear patterns, e.g. the $1$'s are at positions $1,1+2=3, 1+2+3=6, \ldots$; further, those $1$'s are followed by $2,3,4,\ldots$.  Those two patterns are simple to prove, but a general formula for the $n$-th term remains a problem for the future.  We close by mentioning that the Pulsar sequence is, at present, missing from the wonderful repository \url{oeis.org}, the Online Encyclopedia of Integer Sequences.



\begin{thebibliography}{2}
	\bibitem{ctc}
	  The Sudoku Discovery Of The Decade: The Sequel!!!, published 7/17/2025,
      \newblock
    \url{https://www.youtube.com/watch?v=uymIHULB12c}.

	
\end{thebibliography}
\end{document}